\documentclass[11pt,a4paper]{article}
\pdfoutput=1
\long\def\meta#1{\texttt{#1}}
\long\def\drop#1{}

\usepackage[latin1]{inputenc}
\usepackage{amsmath}
\usepackage{amsthm}
\usepackage{amssymb}
\usepackage{latexsym}
\usepackage{graphicx}
\usepackage{mathrsfs}
\usepackage{geompsfi}
\usepackage{psfrag} 
\usepackage{subfig}
\usepackage{bm}

\def\R{\mathbb R}
\def\div{\mathop{\mathrm{div}}}
\let\o\Omega

\let\ra\rightarrow

\let\g\gamma
\newcommand{\hf}{\mathscr{H}^1}
\def\pref#1{(\ref{#1})}
\newtheorem{theorem}{Theorem}

\newtheorem{lemma}[theorem]{Lemma}

\newenvironment{remark}%
  {\par\medbreak\refstepcounter{theorem}%
    \noindent\textbf{Remark~\thetheorem. }}%
  {\par\medskip}

\begin{document}
\title{Non-oriented Solutions of the Eikonal Equation}
\author{Mark A. Peletier and Marco Veneroni}
\date{\today}
\maketitle

\begin{abstract}
We study a new formulation for the eikonal equation $|\nabla u| =1$ on a bounded subset of $\R^2$. Instead of a vector field $\nabla u$, we consider a field $P$ of orthogonal projections on $1$-dimensional subspaces, with $\div P \in L^2$. We prove existence and uniqueness for solutions of the equation $P\div P=0$. We give a geometric description, comparable with the classical case, and we prove that such solutions exist only if the domain is a tubular neighbourhood of a regular closed curve. 
The idea of the proof is to apply a generalized method of characteristics introduced in \cite{JabinOttoPerthame02} to a suitable vector field $m$ satisfying $P=m\otimes m$.

This formulation provides a useful approach to the analysis of stripe patterns. It is specifically suited to systems where the physical properties of the pattern are invariant under rotation over 180 degrees, such as systems of block copolymers or liquid crystals.

\end{abstract}
\smallskip

\noindent\textbf{AMS Cl.} 35L65, 35B65.

\section{Introduction}

\subsection {Stripe patterns and the eikonal equation}

Many pattern-forming systems produce parallel stripes, sometimes straight, sometimes curved. In geology, for instance, `parallel folding' refers to the folding of layers of rock in a manner that preserves the layer thickness but allows for curving of the layers~\cite{BoonBuddHunt07}. In a different context, the convection rolls of the Rayleigh-B\'enard experiment produce striped patterns that may also be either straight or curved (see e.g.~\cite{BodenschatzPeschAhlers00}). The system that suggested the work of this paper is a third example: in~\cite{PeletierVeneroniTR} we investigated striped patterns that arise in the modelling of \emph{block copolymer melts}.  

Block copolymers consist of two covalently bonded, mutually repelling parts (`blocks'). At sufficiently low temperature the repelling forces lead to patterns with a length scale that is related to the length of single polymers. We recently studied the behaviour of an energy that describes such systems, and investigated a limit process in which the stripe width tends to zero~\cite{PeletierVeneroniTR}. In that limit the stripes not only become thin, but also uniform in width, and the stripe pattern comes to resemble the level sets of a solution of the eikonal equation. The rigorous version of this statement, in the form of a Gamma-convergence result, gives rise to a new formulation of the eikonal equation, in which the directionality of the stripes is represented by line fields rather than by vector fields. Before stating this formulation mathematically we first describe it in heuristic terms.

\medskip
The eikonal equation has its origin in models of wave propagation, where the equation describes the position of a wave front at different times $t$. For a homogeneous and isotropic medium, in which the wave velocity is constant, the equation can be written in the form
\begin{equation}
\label{eq:EEvector}
|\nabla u | = 1.
\end{equation}
The wave front at time $t$ is given by the level set $\{x: u(x) = t\}$, and the function $u$ has the interpretation of the time needed for a wave to arrive at the point $x$. 

A feature of the eikonal equation is that the fronts at different times are parallel, provided no singularities occur. In this sense the equation is a natural candidate for the description of other processes that involve parallellism, such as the stripe-forming systems mentioned above. However, a major difference between the stripe-forming systems and the wave-front model is that the wave front has a natural \emph{directionality} associated with it: of the two directions normal to a front, one is `forward in time' and the other `backward'. This distinction also is visible in the notion of viscosity solution for~\pref{eq:EEvector} (see e.g.~\cite{CrandallEvansLions84}).

\begin{figure}[ht]
\centering
\noindent
\includegraphics[height=2cm]{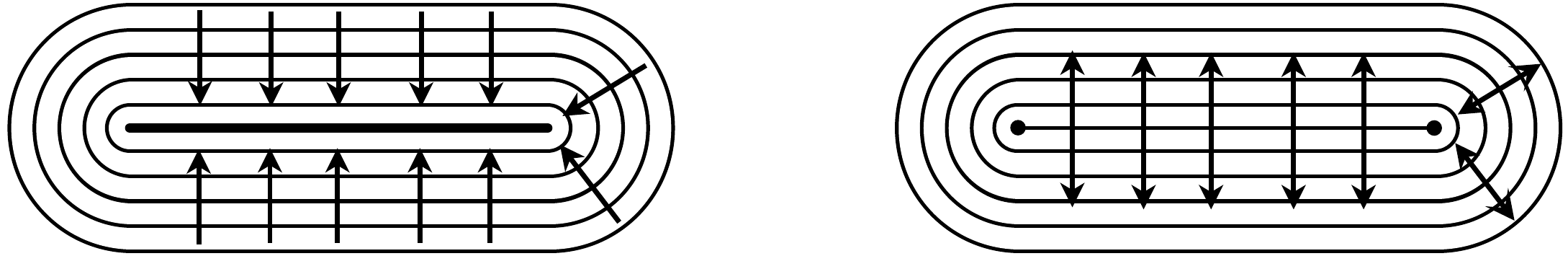}
\caption{Stripe patterns can be represented by vectors (left) or by unoriented line fields (right). Both representations have a vortex singularity at the two ends; but the vector representation also contains a jump singularity along the connecting line. Note that the regularity restrictions of this paper exclude both types of singularity, however.}
\label{fig:stadium}
\end{figure}

The stripe patterns, on the other hand, have no inherent distinction between the two normal directions. As a consequence a vector representation of a stripe pattern may have singularities that have no physical counterpart. Figure~\ref{fig:stadium} (left) shows an example of this.

In our work on block copolymers that was mentioned above~\cite{PeletierVeneroniTR}, the limiting stripe problem is formulated in terms of  \emph{unsigned vector fields}, or \emph{line fields}, which capture direction only up to a sign (Figure~\ref{fig:stadium} (right)). These stripe patterns are required to satisfy what is essentially \emph{the eikonal equation for line fields} (see~(\ref{pb:main}) below). In this paper we investigate this new formulation of the eikonal equation, study the properties of its solutions, and compare it to various existing formulations. We only consider the case of two space dimensions. 

\subsection{The eikonal equation in terms of line fields}
We now turn to the work of this paper. A natural mathematical object for the representation of line fields is a \emph{projection}. For the purposes of this paper, we define a {projection} to be 
a matrix $P$ that can be written in terms of a unit vector $m$ as $P=m\otimes m$. 
Such a projection matrix has a range and a kernel that are both one-dimensional, and if necessary one can identify a projection $P$ with its range, i.e. with the one-dimensional subspace of $\R^2$ onto which it projects. Note that the independence of the sign of $m$, that is the unsigned nature of a projection, can be directly recognized in the formula $P=m\otimes m$.

We define $\div P$ as the vector-valued function whose $i$-th component is given by $(\div P)_i:= \sum_{j=1}^2 \partial_{x_j}P_{ij}$. We consider the following problem. Let $\Omega$ be an open subset of $\R^2$.
 Find $P\in L^\infty(\Omega;\R^{2\times 2})$ such that
\begin{subequations}
\label{pb:main}
\begin{eqnarray}
         P^2 = P && \mbox{a.e. in } \o,\label{prpro}\\
         \mbox{rank}(P)=1 && \mbox{a.e. in } \o,\label{prrank}\\
         P \mbox{ is symmetric} && \mbox{a.e. in } \o,\label{prsymm}\\
         \div P \in L^2(\R^2;\R^2) && (\mbox{extended to $0$ outside  }\o),\label{divpint}\\
         P\, \div\! P = 0  && \mbox{a.e. in } \o.\label{pdivpzero}
\end{eqnarray}
\end{subequations}
The first three equations encode the property that $P(x)$ is a projection, in the sense above, at almost every $x$. Equation~\pref{divpint} is both a regularity requirement and a boundary condition, and the choice of the exponent $2$ has a precise explanation, as we show below. 

Given the regularity provided by~\pref{divpint}, the final condition~\pref{pdivpzero} is the eikonal equation itself, as a calculation for a smooth unit-length vector field $m(x)$ shows: 
\begin{equation}
\label{comp:PdivP}
0 = P\div P = m (m\cdot (m\div m + \nabla m \cdot m))
= m\div m + m(m\cdot \nabla m \cdot m)
= m\div m,
\end{equation}
where the final equality follows from differentiating the identity $|m|^2 =1 $. A solution vector field $m$ therefore is divergence-free, implying that its rotation over 90 degrees is a gradient $\nabla u$; from $|m|=1$ it follows that $|\nabla u|=1$ (see Section \ref{sec:potential}). This little calculation also shows that the 
interpretation of $m$ in $P=m\otimes m$ is that of the stripe direction; $P$ projects along the normal onto the tangent to a stripe.

The sense of property (\ref{divpint}) is that the divergence of $P$ (extended to $0$ outside $\o$), in the sense of distributions in $\R^2$, is an $L^2(\R^2)$ function, i.e. there exists $C>0$ such that for any test function $\varphi \in C^{\infty}_c(\R^2,\R^2)$
\begin{equation}\label{distr}
    \left|\int_{\R^2}P(x):\nabla \varphi(x)\,dx\right| \leq C{\|\varphi\|}_{L^2(\R^2)}.
\end{equation}
Since the trace $Pn$ satisfies the equality 
\begin{equation}
\label{calc:bdry}
         -\int_\o P:\nabla\varphi\, dx = \int_\o \div P \cdot \varphi\, dx - \int_{\partial\o}\!(P n)\cdot \varphi\, dS,
\end{equation}
condition~\pref{divpint} implies 
\begin{equation}\label{ptgbound}
        Pn = 0\quad \mbox{in the sense of traces on }\partial\o.
\end{equation}
In order to see this, note that by \pref{distr} and \pref{calc:bdry} there exists $C>0$ such that
$$ \left|\int_\o \div P \cdot \varphi\, dx - \int_{\partial\o}\!(P n)\cdot \varphi\, dS \right| \leq C{\|\varphi\|}_{L^2(\R^2)},  \quad \forall\, \varphi \in C^{\infty}_c(\R^2,\R^2).$$
Choosing a sequence $\{\varphi_k\}$ such that $\varphi_k \ra 0$ in $L^2(\R^2)$ and $\varphi_k \ra Pn$ on $L^2(\partial \o)$, we conclude \pref{ptgbound}.
\medskip

The exponent $2$ in~\pref{divpint} is critical in the following sense. Obvious possibilities for singularities in a line field are jump discontinuities (`grain boundaries') and target patterns (see Figure~\ref{fig:types_of_variation}). 
\def\hht{2cm}%
\def\spacing{\hskip1cm}%
\begin{figure}[ht]
\centering
\subfloat[caption][\centering grain boundary]{\includegraphics[height=\hht,clip=on]{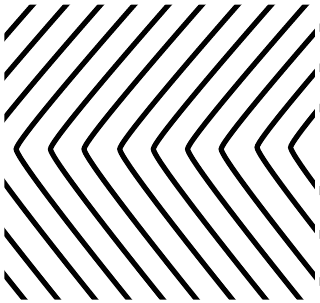}\label{subfig:gb}}
\spacing
\subfloat[caption][\centering target and U-turn patterns]{\includegraphics[height=\hht,clip=on]{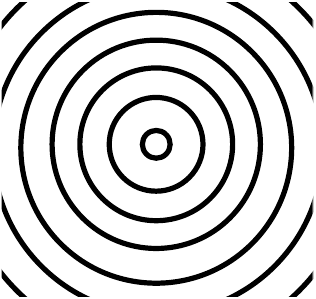}\ \ \includegraphics[height=\hht,clip=on]{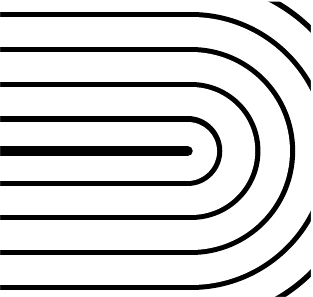}\label{subfig:target}}
\spacing
\subfloat[caption][\centering smooth directional variation]{\includegraphics[height=\hht,clip=on]{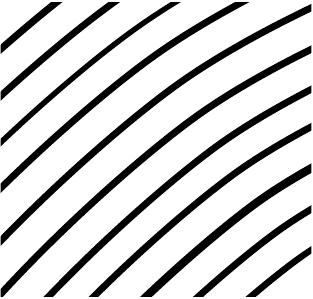}}
\caption{Canonical types of stripe variation in two dimensions. Types~\protect\subref{subfig:gb} and~\protect\subref{subfig:target} are excluded by~\pref{divpint}.
}
\label{fig:types_of_variation}
\end{figure}%
At a grain boundary the jump in $P$ causes $\div P$ to have a line singularity, comparable to the one-dimensional Hausdorff measure; condition~\pref{divpint} clearly excludes that possibility. For a target pattern the curvature $\kappa$ of the stripes scales as $1/r$, where $r$ is the distance to the center; then  $\int \kappa^p$ is locally finite for $p<2$, and diverges logarithmically for $p=2$. The cases $p<2$ and $p\geq 2$ therefore distinguish between whether target patterns are admissible ($p<2$) or not.

\subsection{Results}
The first main contribution of this paper is to show that problem~\pref{pb:main} indeed resembles the eikonal equation. In the heuristic discussion above, we argued that the field $P$ has to be locally parallel, to be parallel to the boundary of $\Omega$, and to avoid line and vortex singularities. 
The two first statements are formalized in the following theorem.
\begin{theorem}
\label{th:ee0}
Let $\Omega$ be an open, bounded, and connected subset of $\R^2$ with $C^2$ boundary, and let $P$ be a solution of~\pref{pb:main}.
\begin{enumerate}
\item Let $x_0$ be a Lebesgue point of $P$ in $\Omega$, let $x\in \Omega$, and let $L$ be the line segment connecting $x_0$ with $x$. Assume that $L\subset \Omega$. If $P(x_0)\cdot (x-x_0) =0$, then $P(y)=P(x_0)$ for $\hf$-almost every $y\in L$.
\label{th:ee0:orth}
\item $P\cdot n=0$ a.e. on $\partial \Omega$.
\label{th:ee0:bdry}
\end{enumerate}
These two statements are meaningful since 
\begin{enumerate}
\setcounter{enumi}2
\item $P\in H^1(\Omega;\R^{2\times2})$.
\label{th:ee0:regularity}
\end{enumerate}
\end{theorem}
We recall that $x$ is called a Lebesgue point of $m$ if 
$$ \lim_{r \ra 0^+} \frac1{r^2}\int_{B_r(x)}|m(x)-m(y)|dy =0.$$
A standard result yields that a.e. $x\in \o$ is a Lebesgue point for $m$ (see e.g.~\cite[Section 1.7]{EvansGariepy92}).

\bigskip

The second main result is to show that the restrictions on $P$ are so rigid that the mere existence of a solution provides a strong characterization of the \emph{geometry of the domain}~$\Omega$:

\begin{theorem}
\label{th:ee}
Let $\Omega$ be an open, bounded, and connected subset of $\R^2$ with $C^2$ boundary. Then there exists a solution of~\pref{pb:main} if and only if $\Omega$ is a tubular domain. In that case the solution is unique.
\end{theorem}

A \emph{tubular domain} is a domain in $\R^2$ that can be written as
\[
\Omega = \Gamma + B(0,\delta),
\]
where $\Gamma$ is a closed curve in $\R^2$ with continuous and bounded curvature $\kappa$, $0<\delta< \|\kappa\|^{-1}_\infty$, and $B(0,\delta)$ is the open ball of center 0 and radius $\delta$.

\medskip

\begin{figure}[ht]
\indent\centering
\psfig{figure=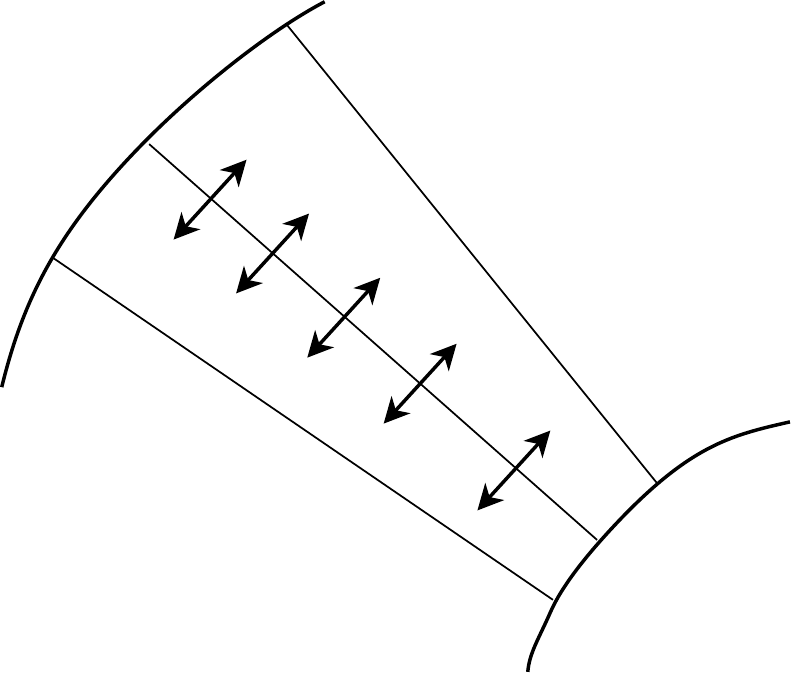,height=4cm}
\caption{If tangent directions propagate normal to themselves (Theorem \ref{th:ee0}.\ref{th:ee0:orth}), and if in addition the boundary is a tangent direction (\ref{ptgbound}), then the domain is tubular (Theorem~\ref{th:ee}).}
\label{fig:ee}
\end{figure}
The reason why Theorem~\ref{th:ee} is true can heuristically be recognized in a simple picture. Figure~\ref{fig:ee} shows two sections of $\partial \Omega$ with a normal line that connects them. By the first assertion of Theorem~\ref{th:ee0}, the stripe tangents are orthogonal to this normal line; by the second, this normal line is orthogonal to the two boundary segments, implying that the two segments have the same tangent. Therefore the length of the connecting normal line is constant, and as it moves it sweeps out a full tubular neighbourhood.

\subsection{Discussion}

The work of this paper represents a first step in the analysis of this projection-valued eikonal  equation. While the main results are still lacking in various ways---which we discuss in more detail below---the main point of this paper is to show that this projection-valued formulation is a useful alternative to the usual vector-based formulation. 

To start with, our Theorems~\ref{th:ee0} and~\ref{th:ee} show that solutions of~\pref{pb:main} behave much like we expect from the eikonal equation, in the sense that directional information is preserved in the normal direction. Theorem~\ref{th:ee} makes this property even more explicit, by showing that a full tube, or bunch, of parallel `stripes' can be identified. 

However, it is the differences with the vector-valued eikonal equation that are the most interesting. Figure~\ref{fig:stadium} shows how this formulation can be a better representation of the physical reality than the vector-based form. On the left, the vector field has a jump discontinuity along the center line, while on the right the projection is continuous along that line. Depending on the underlying model, this singularity may have a physical counterpart, or may be a spurious consequence of the vector-based description. For the wave-propagation model the singularity is very real; for striped-pattern systems it typically is not. A projection-valued formulation therefore provides an alternative to the Riemann-surface approach that is sometimes used~\cite{ErcolaniIndikNewellPassot03}. For this distinction to have any consequence, however, solutions with less regularity than the $\div P\in L^2$ of this paper are to be considered. 

\medskip
We now comment in more detail on our method of proof. The proof of the properties that we give in this paper relies on a reduction of the projection-valued formulation to a vector-based formulation. This reduction is achieved by the Ball-Zarnescu lemma (Lemma~\ref{bz}), which requires $\div P\in L^2$; for less regularity the existence of a lifting may not hold, as the example of the U-turn pattern (Figure~\ref{subfig:target}) shows. 

The dependence of the proof on a vector-based representation is awkward in various ways. To start with, the condition $\div P\in L^2$ required for the lifting is much stronger than the conditions~(\ref{normone}--\ref{kine}) that Jabin, Otto, and Perthame require for their results~\cite{JabinOttoPerthame02}. It also has the effect of excluding all singularities, as we already remarked. It would be interesting to prove properties such as those of Theorems~\ref{th:ee0} and~\ref{th:ee} by methods that do not rely on this lifting. 

We would hope that such an intrinsic projection-based proof could also be generalized to the study of target patterns and U-turns, and eventually of grain boundaries. These will require increasingly weak regularity requirements: target patterns may exist for $\div P\in L^p$ with $p<2$, and for a line discontinuity, such as a grain boundary, $\div P$ will be a measure.

\subsection{Overview of the paper}

The structure of the rest of the paper is as follows. We start, in Section~\ref{sec:regularity}, by showing that $P\in H^1$ (part~\ref{th:ee0:regularity} of Theorem~\ref{th:ee0}). In Section~\ref{sec:lifting} we construct a lifting, i.e. a vector field $m$ such that $P=m\otimes m$; there we also show part~\ref{th:ee0:bdry} of Theorem~\ref{th:ee0}. In Section~\ref{sec:properties} we show that the lifting $m$ satisfies the conditions of~\cite{JabinOttoPerthame02}, and prove the final part of Theorem~\ref{th:ee0}. In Section~\ref{sec:th2} we combine all the earlier results to prove Theorem~\ref{th:ee}.

\section{Regularity of $P$}
\label{sec:regularity}

The regularity statement~\ref{th:ee0:regularity} of Theorem~\ref{th:ee0} follows from simple manipulation. 
By (\ref{prpro}--\ref{prsymm}) $P$ is a symmetric matrix which we can write as
$$P=\left(
\begin{array}{cc}
    a & b \\
    b & c
\end{array}
\right),
$$
and whose elements satisfy
\begin{gather}\label{abcbound}
        a,b,c\in L^\infty(\Omega), \quad 0\leq a,b,c \leq 1,\\
\label{elementary}
    b^2=ac, \quad a+c =1.
\end{gather}
Denote by $a_i:=\partial_{x_i}a$ (idem for $b_i,c_i$) for $i=1,2$. By (\ref{divpint}) there exist $f,g\in L^2(\o)$ such that
\begin{equation}\label{abcl}
    \left\{
        \begin{array}{l}
            a_1 + b_2 = f\\
            b_1 + c_2 = g.
        \end{array}
      \right.
\end{equation}
From~\pref{elementary} we deduce by differentiation
\[
bb_1 = \left(\frac12 - a\right)a_1,
\]
and by using this and  (\ref{abcbound}--\ref{abcl}) in \pref{pdivpzero} we have 
\begin{eqnarray*}
    0   &=& a(a_1 + b_2) + b(b_1+c_2)\\
        &=& a(a_1 + b_2) + \left(\frac12 -a\right)a_1 +bc_2\\
        &=& \frac12 a_1 + ab_2 -ba_2    \\
        &=& \frac12 a_1 + a(f -a_1) -ba_2.
\end{eqnarray*}
Similarly
\begin{eqnarray*}
    0   &=& b(a_1 + b_2) + c(b_1+c_2)\\
        &=& ba_1 +\left(\frac12 -c\right)c_2 + c(b_1+c_2)\\
        &=& ba_1 -\frac12 a_2 + (1-a)b_1\\
        &=& ba_1 -\frac12 a_2 +  (1-a)(g + a_2).
\end{eqnarray*}
Therefore $a_1,a_2$ solve the linear system
$$
\left(
    \begin{array}{cc}
        \frac12 -a & -b\\
        b       &   \frac12 -a
    \end{array}
\right)
\left(
    \begin{array}{c}
        a_1\\
        a_2
    \end{array}
\right)
=
\left(
    \begin{array}{c}
        -a f\\
        (a-1)g
    \end{array}
\right).
$$
Since the determinant of the coefficient matrix is $(\frac12 -a)^2 + b^2=\frac14>0$, the system is nondegenerate, and it follows from $a,b,c\in L^{\infty}(\o)$ and $f,g\in L^2(\o)$ that $a_1,a_2\in L^2(\o).$ By the relations (\ref{elementary}) and (\ref{abcl}) we find $c_1,c_2\in L^2(\o)$ and $b_1,b_2\in L^2(\o)$. Therefore
$$ P_{ij}\in H^1(\o),\quad i,j=1,2.$$ 
\qed

\section{Lifting to a vector representation}
\label{sec:lifting}

The remaining parts of both theorems are proved by using a vectorial representation of~$P$. We say that $P$ is \emph{orientable} if there exists a vector field $m\in H^1(\o;S^1)$ such that 
\begin{equation}\label{orientable}
    P=m \otimes m\quad \mbox{a.e. on } \o.
\end{equation}
In this case $m$ is called a \emph{lifting} of $P$. Similarly, the trace of $P$ on the boundary is defined to be {orientable} if there exists a vector field  $\overline{m}\in H^{1/2}(\partial\o;S^1)$, such that
\begin{equation}\label{trorient}
        \mbox{Tr}(P)_{|\partial\o}=\overline{m}\otimes \overline{m} \quad \text{a.e. on }\partial \Omega.
\end{equation}

The following lemma by Ball and Zarnescu establishes a link between orientability in the bulk and on the boundary. We formulate it in the language of this paper. 
\begin{lemma}[\cite{BallZarnescu07TA}]\label{bz}
Let $\o$ be an open, bounded,  connected subset of $\R^2$, with $C^2$ boundary. If $P$ satisfies~\pref{pb:main} and has the additional regularity $P\in  H^1(\o;\R^{2 \times 2})$, then $P$  is orientable if and only if its trace on $\partial\o$ is orientable.
\end{lemma}

Note that the boundary $\partial\Omega$ has only a finite number $N\geq 1$ of connected components. Indeed, let $\Lambda$ be a set of indexes, and let $\g_\lambda$, for $\lambda\in \Lambda$, be a connected component of $\partial \o$, which disconnects $\R^2$ into a bounded set $B_\lambda$ and an unbounded set $U_\lambda$. Since $\o$ is connected, 

\begin{itemize}
	\item there is a unique index $\bar{\lambda}$ such that $B_{\bar{\lambda}}\cap \o \neq \emptyset$ and 
	\item $B_{\lambda_1}\cap B_{\lambda_2}=\emptyset$ if $\lambda_1\neq \lambda_2$ and $\lambda_1,\lambda_2\neq \bar{\lambda}$. 
\end{itemize}
Since $\partial\o$ is $C^2$, its curvature is bounded by some constant $\delta$ and therefore every $B_\lambda$ must contain at least one ball of radius $1/\delta$. Since $\partial \o$ is bounded, $\bigcup_{\lambda\in \Lambda, \lambda \neq \bar{\lambda}} B_\lambda$ is bounded and only a finite number of such balls can fit into it without overlapping. We conclude that $\# \Lambda <\infty.$

We now construct a lifting $\overline m$ of $P$ on the boundary. Let $\partial\o=\partial\o^0\cup\dots\cup\partial\o^N$ be a decomposition of the boundary into connected components, and let $\alpha_j:[0,L_j]\ra \partial\o^j$, be $C^2$ arclength parameterizations of $\partial\o^j$. Then, owing to (\ref{ptgbound}), the vectors $\overline m^j$ defined by $\overline{m}^j(\alpha_j(s)):=\alpha_j'(s)$, $j=0\ldots N$, satisfy (\ref{trorient}). Therefore, by Lemma \ref{bz} we obtain that $P$ is orientable.

\medskip
Thus we have proved 
\begin{lemma}
There exists a vector field $m\in H^1(\o,\R^2)$ such that $ P = m \otimes m$.
\end{lemma}
This proves part~\ref{th:ee0:regularity} of Theorem~\ref{th:ee0}. Part~\ref{th:ee0:bdry} then follows from the calculation~\pref{calc:bdry}. Note that $Pn=0$ on the boundary implies
\begin{equation}
\label{eq:mbdry}
m \cdot n =0 \qquad\text{a.e. on }\partial \Omega.
\end{equation}

\section{Properties of $m$}
\label{sec:properties}

We use some concepts and notation from~\cite{JabinOttoPerthame02}. 
We define the functions $\chi$ and ${\bm \chi}$ by
\begin{equation}\label{defchi}
   \chi(x,\xi):= {\bm \chi}(m(x),\xi) := \left\{ \begin{array}{ll}
       1 & \mbox{if } m(x)\cdot \xi> 0\\
       0 & \mbox{if } m(x)\cdot \xi\leq 0.
   \end{array}
\right. 
\end{equation}
\begin{lemma}\label{ptom}
The vector field $m$ satisfies
            \begin{eqnarray}
                |m(x)| &=\ 1& \mbox{ for a.e. $x$ in } \o,\label{normone}\\
                \div \,m &=\ 0& \mbox{ distributionally in } \R^2,\label{divzero}\\
                \xi \cdot \nabla \chi(\cdot,\xi) &=\ 0& \mbox{ distributionally in } \o \mbox{ for all $\xi$ in } S^1.\label{kine}
            \end{eqnarray}
\end{lemma}

\begin{proof}
Property (\ref{normone}) follows from remarking that
\[
|m|^4 = m\cdot P \cdot m = m\cdot P^2 \cdot m = |m|^6.
\]
The computation~\pref{comp:PdivP} yields~\pref{divzero}. The proof of (\ref{kine}) is an application of some ideas introduced in \cite{DeSimoneKohnMuellerOtto01} and further developed in \cite{JabinOttoPerthame02}. We use the notion of entropy introduced in \cite[Definition 2.1]{DeSimoneKohnMuellerOtto01}. A vector field $\Phi \in C^{\infty}_0(\R^2,\R^2)$ is called an \emph{entropy} if
\begin{equation}\label{entropy}
    z \cdot \nabla\Phi z^\perp =0\quad \mbox{for all }z\in \R^2,\quad \mbox{and}\quad \Phi(0)=0,\ \nabla\Phi(0)=0.
\end{equation}
In \cite{DeSimoneKohnMuellerOtto01} and \cite{JabinOttoPerthame02} it was shown that for every entropy $\Phi$ there exists a vector field $\Psi\in C^{\infty}_0(\R^2,\R^2)$ and a function $\alpha\in C^{\infty}_0(\R^2,\R)$ such that for any $m\in H^1(\o,\R^2)$
\begin{equation*}
    \div(\Phi(m))=\Psi(m)\cdot \nabla(1-|m|^2) + \alpha(m)\,\div m\quad \mbox{a.e. in }\o.
\end{equation*}
As proved in \cite[Lemma 5]{DeSimoneKohnMuellerOtto01}, for any fixed $\xi\in S^1$, the function
\begin{equation}\label{appr}
    \R^2 \ni z \mapsto \bm\chi(z,\xi)\xi
\end{equation}
(where $\bm\chi$ is defined in (\ref{defchi})) is the pointwise limit of a sequence $\{\Phi_n\}_{n \in \mathbb{N}}$ of entropies in the sense of (\ref{entropy}). Properties (\ref{normone}--\ref{divzero}) and the approximation (\ref{appr}) then yield equation~(\ref{kine}).
\end{proof}

The properties of $m$ stated in Lemma \ref{ptom} allow us to apply Proposition 3.2 in \cite{JabinOttoPerthame02}, which reads 
\begin{lemma}[\cite{JabinOttoPerthame02}]\label{propjop} Let $m$ satisfy (\ref{normone}--\ref{kine}). 
        Let $x_0\in \overline\Omega$ 
        be a Lebesgue point for $m$ in $\Omega$, and let $L\subset \overline\Omega$ be a straight line segment containing $x_0$. Then
        $$ m(x)\cdot m^\perp(x_0)=0\quad \mbox{for $\hf$-a.e. } x\in L. $$
\end{lemma}
In terms of $P=m\otimes m$ this statement reduces to part~\ref{th:ee0:orth} of Theorem~\ref{th:ee0}. This concludes the proof of Theorem~\ref{th:ee0}.

\section{Proof of Theorem~\ref{th:ee}}
\label{sec:th2}

The statement of equivalence in Theorem~\ref{th:ee} contains a trivial and a non-trivial part. The non-trivial part is to show that existence of a solution $P$ implies that $\Omega$ is tubular. The trivial part is to construct a solution $P$, if one assumes that $\Omega$ is tubular. We prove the non-trivial part first, since the calculations will be useful in the second part. Part~\ref{th:ee0:orth} of Theorem~\ref{th:ee0} will appear as an intermediate result in the proof of Theorem~\ref{th:ee}.

\subsection{A first characterization of $\Omega$}

Let $\mathcal{A}$ be the class of sets $\o$ such that
\begin{quote}
    \sl
    $\o \subset \R^2$ is open, bounded, connected, $C^2$, and $\exists\,y,z \in \partial\o$ such that the normal lines
    issued from $y,z$ are different and intersect in $\o$ before crossing~$\partial\o$.
\end{quote}
\smallskip

Then Theorem 1.2 in \cite{JabinOttoPerthame02} states that if $\o\in \mathcal{A}$ and $m$ satisfies (\ref{normone}), (\ref{divzero}), and (\ref{kine}), then $\o$ is a disk and $m$ is a vortex, i.e. there exist $x_0\in \o$ and $\alpha=\pm 1$ such that 
\begin{equation}
\label{def:vortex}
m(x)=\alpha\frac{(x-x_0)^\perp}{|x-x_0|}\quad \mbox{for a.e. }x\in\o. 
\end{equation}
(Here and below, $v^\perp$ is the counterclockwise rotation by 90 degrees of the vector $v$.)
A direct computation shows that if $P=m\otimes m$ with $m$ given by~\pref{def:vortex} then $\div P \notin L^2(\o)$.
It follows that $\o\notin \mathcal{A}$, and therefore either the normal lines issued from any two points on the boundary are identical or they have no intersection in $\o$. 

\subsection{Construction of a potential}\label{sec:potential}
Let $\g:[0,L]\ra \partial\o^0$ be a $C^2$ arclength parameterization of $\partial\o^0$ and let $n$ be the exterior unit normal vector on~$\partial\o$. By property (\ref{divzero}) the vector field $m^\perp$ satisfies
$$ 
\mbox{curl }\, m^\perp =0,\ \mbox{ distributionally in } \R^2.
$$
Since the vector field $m^\perp$ is orthogonal to the boundary (see~\pref{eq:mbdry}), by Green's theorem the integral of $m^\perp$ along any one-dimensional closed curve in $\Omega$ is zero. Therefore
there exists a potential $\phi\in H^2(\o)$ which satisfies
\begin{equation}\label{phiedm}
        \nabla \phi(x) = m^\perp(x).
\end{equation}
We can assume that $\phi$ is continuous on $\overline \Omega$. 
From~\pref{eq:mbdry} we deduce
\begin{equation}\label{perp}
    \nabla\phi \cdot n^\perp = -m\cdot n = 0\quad \mbox{a.e. on }\partial\o,
\end{equation}
therefore $\phi$ is constant on every connected component of $\partial\o$: i.e. there exist $c_j\in \R$, $j=0,\ldots,N$ such that
\begin{equation}
\label{phiconst}
\phi(x)=c_j\quad \forall\, x\in \partial\o^j,\ j=1,\ldots,N.
\end{equation}
It is not restrictive to assume $c_j\geq 0,$ $c_0=0$. Moreover for $x\in\partial\o^0$ and $T>0$ we define the segments
$$
L_{x,T}:= \{x-tn(x):t\in[0,T]\},
$$
which intersect $\partial\o$ orthogonally by (\ref{perp}) (wherever they intersect), and the function $T(x):\partial\o^0\ra \R$,
$$ 
T(x):=\max\{T>0: L_{x,T}\subset \overline{\o}\}.
$$
Note that it follows from this definition that 
\begin{equation}
\label{prop:TonpOm}
x-T(x)n(x)\in \partial\Omega.
\end{equation}

\subsection{$\Omega$ is tubular}

By \cite[Theorem 2, section 5.3]{EvansGariepy92}, there exists a set $\mathcal{N}\subset \partial\o$ such that $\hf(\mathcal{N})=0$ and every $x\in \partial\o \backslash \mathcal{N}$ is a Lebesgue point for $m$ with respect to the two-dimensional Lebesgue measure. 

We now proceed in five steps.

\textit{i)}\ 
As a direct consequence of Lemma \ref{propjop},
$$
\forall\,x\in\partial\o^0\backslash \mathcal{N}, \mbox{ for $\hf$-a.e. }y\in L_{x,T(x)},\ m(y)=\pm m(x);
$$

\textit{ii)}\ We claim, in addition, that
$$\forall\,x\in\partial\o^0\backslash \mathcal{N}, \mbox{ for $\hf$-a.e. }y\in L_{x,T(x)},\ m(y)= m(x);$$
otherwise, following the arguments in the proof of \cite[Theorem 1.1]{JabinOttoPerthame02}, it can be seen that a point where $m(y)$ jumps from $m(x)$ to $-m(x)$ would be a vortex point for $m$; as we remarked above, that would imply that $\div P\not\in L^2$, a contradiction. 
\medskip

\textit{iii)}\ By (\ref{phiedm}) and step \textit{ii)}, for each $x_0\in \partial \Omega^0\setminus\mathcal N$, $\phi$ is linear on $L_{x_0,T(x_0)}$ with slope $1$.

\textit{iv)} By~\pref{prop:TonpOm} and step \textit{iii)}, for $\hf$-a.e. $x\in \partial\o^0$ there exists a $j\in \{1,\ldots,N\}$ such that
$$ 
c_j - 0 = \phi\big(x-T(x)\,n(x)\big)-\phi(x)=T(x).
$$
Therefore there exists a step function $\bar{T}:\partial\o^0 \ra \{c_1,\ldots,c_N\}$ such that
$$ T = \bar{T},\quad \hf\mbox{-a.e. on }\partial\o^0.$$

\textit{v)}\ Let $\bar{x}$ be a discontinuity point  for $T$.

\begin{figure}[ht]
\centering \noindent
\psfig{figure=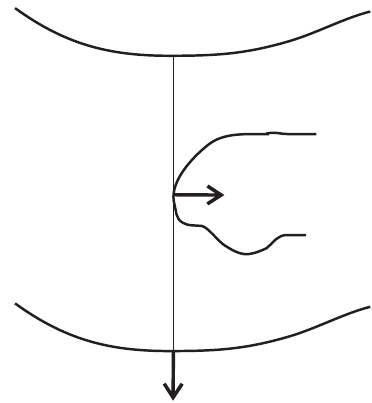,height=3.5cm}\caption{For example, $T$ in $\bar{x}$ jumps between the values $c_1$ and $c_2$.}
\label{fig:jump}
\end{figure}
We claim that there exists $\hat{x}\in L_{\bar{x},T(\bar{x})}\cap \partial\o$ such that $n(\bar{x})\cdot n(\hat{x})=0;$
otherwise, the regularity of $\partial\o$ and the continuity of $n$ would make  $T$ continuous in $\bar{x}$. This contradicts $\o\notin \mathcal{A}$, and we conclude that $T$ is continuous and $T(x)\equiv c$. Finally, the image of $\partial\o^0$ under the continuous function $x\mapsto x+T(x)\,n(x)$ is connected, and therefore $\partial \o$ consists of two connected components at constant distance $c$: $\partial\o=\partial\o^0\cup\partial\o^1$ and $\o=\Gamma + B(0,\delta)$ for $\delta=c/2$ and  $\Gamma:=\{ x+\frac c2 n(x), x\in\partial\o^0\}$. This concludes the proof of one half of Theorem~\ref{th:ee}.
\rightline{$\Box$}

\subsection{Proof of the converse}
Let $\o$ be a tubular domain of width $\delta$ and let $\partial\o = \partial\o^0 \cup \partial\o^1$ be the decomposition of the boundary into connected components. Let $\g:[0,L]\ra \partial\o^0$ be a $C^2$ arclength parameterization of $\partial\o^0$, then
\begin{align*}
    \Psi &: [0,L]\times (0,2\delta)\ra \o\\
    \Psi(s,t)&:= \gamma(s) -t n(s)
\end{align*}
is a $C^1$ parameterization of $\o$. Then the line field $P$ defined as
\begin{equation}\label{formula}
    P(\Psi(s,t)):=\g'(s)\otimes \g'(s)
\end{equation}
is a solution and satisfies $P\in C^1(\o,\mathbb{M}^{2 \times 2})$. Since by step \textit{ii)} above $m$ is uniquely determined a.e. on $\o$ by its trace on  $\partial\o$, i.e. by the tangents to $\partial\o$, we conclude that $P$ is uniquely determined a.e. on $\o$ and (\ref{formula}) is the unique solution.
\rightline{$\Box$}
\smallskip

\bibliographystyle{plain}
\bibliography{eikon}

\begin{thebibliography}{1}

\bibitem{BallZarnescu07TA}
J.~M. Ball and A.~Zarnescu.
\newblock Orientable and non-orientable director fields for liquid crystals.
\newblock To appear.

\bibitem{BodenschatzPeschAhlers00}
Eberhard Bodenschatz, Werner Pesch, and Guenter Ahlers.
\newblock Recent developments in {R}ayleigh-{B}\'enard convection.
\newblock In {\em Annual review of fluid mechanics, {V}ol. 32}, volume~32 of
  {\em Annu. Rev. Fluid Mech.}, pages 709--778. Annual Reviews, Palo Alto, CA,
  2000.

\bibitem{BoonBuddHunt07}
J.~A. Boon, C.~J. Budd, and G.~W. Hunt.
\newblock Level set methods for the displacement of layered materials.
\newblock {\em Proc. R. Soc. Lond. Ser. A Math. Phys. Eng. Sci.},
  463(2082):1447--1466, 2007.

\bibitem{CrandallEvansLions84}
M.~G. Crandall, L.~C. Evans, and P.-L. Lions.
\newblock Some properties of viscosity solutions of {H}amilton-{J}acobi
  equations.
\newblock {\em Trans. Amer. Math. Soc.}, 282(2):487--502, 1984.

\bibitem{DeSimoneKohnMuellerOtto01}
Antonio DeSimone, Stefan M{\"u}ller, Robert~V. Kohn, and Felix Otto.
\newblock A compactness result in the gradient theory of phase transitions.
\newblock {\em Proc. Roy. Soc. Edinburgh Sect. A}, 131(4):833--844, 2001.

\bibitem{ErcolaniIndikNewellPassot03}
N.~Ercolani, R.~Indik, A.~C. Newell, and T.~Passot.
\newblock Global description of patterns far from onset: a case study.
\newblock In {\em Nonlinear dynamics ({C}anberra, 2002)}, volume~1 of {\em
  World Sci. Lect. Notes Complex Syst.}, pages 411--435. World Sci. Publ.,
  River Edge, NJ, 2003.

\bibitem{EvansGariepy92}
L.~C. Evans and R.~F. Gariepy.
\newblock {\em Measure Theory and Fine Properties of Functions}.
\newblock Studies in Advanced Mathematics. CRC Press, 1992.

\bibitem{JabinOttoPerthame02}
Pierre-Emmanuel Jabin, Felix Otto, and Beno{\^{\i}}t Perthame.
\newblock Line-energy {G}inzburg-{L}andau models: zero-energy states.
\newblock {\em Ann. Sc. Norm. Super. Pisa Cl. Sci. (5)}, 1(1):187--202, 2002.

\bibitem{PeletierVeneroniTR}
M.~A. Peletier and M.~Veneroni.
\newblock Stripe patterns in a model for block copolymers.
\newblock To appear.

\end{thebibliography}

\pagebreak

\noindent Mark A. Peletier\smallskip

\noindent Dept. of Mathematics and Computer Science\\
Technische Universiteit Eindhoven\\
PO Box 513\\
5600 MB  Eindhoven\\
The Netherlands\\
e-mail: \meta{m.a.peletier@tue.nl}
\vspace{1cm}

\noindent Marco Veneroni\smallskip

\noindent Dept. of Mathematics and Computer Science\\
Technische Universiteit Eindhoven\\
PO Box 513\\
5600 MB  Eindhoven\\
The Netherlands\\
e-mail: \meta{m.veneroni@tue.nl}

\end{document}